\documentclass[11pt,reqno]{amsart}
\usepackage{amsmath,amscd}
\usepackage{amsthm}
\usepackage{graphicx}
\usepackage{bbm,amssymb}
\usepackage{bbold}
\usepackage{graphicx}
\usepackage[colorlinks]{hyperref}

\newcommand{\floor}[1]{\left\lfloor {#1} \right\rfloor}
\newcommand{\ceiling}[1]{\left\lceil {#1} \right\rceil}

\newcommand{\modulo}[1]{\quad (\mbox{mod }{#1})}

\newtheorem{theorem}{Theorem}
\newtheorem{prop}[theorem]{Proposition}
\newtheorem{lemma}[theorem]{Lemma}
\newtheorem{corollary}[theorem]{Corollary}

\theoremstyle{remark}

\theoremstyle{definition}
\newtheorem*{definition}{Definition}

\def\Leb{\mathcal{L}}
\def\mod{\mbox{mod }}
\def\rot{\rho}

\DeclareSymbolFont{AMSb}{U}{msb}{m}{n}
\DeclareMathSymbol{\C}{\mathbin}{AMSb}{"43} 
\DeclareMathSymbol{\EE}{\mathbin}{AMSb}{"45} 
\DeclareMathSymbol{\N}{\mathbin}{AMSb}{"4E} 
\DeclareMathSymbol{\PP}{\mathbin}{AMSb}{"50} 
\DeclareMathSymbol{\Q}{\mathbin}{AMSb}{"51} 
\DeclareMathSymbol{\R}{\mathbin}{AMSb}{"52} 
\DeclareMathSymbol{\Z}{\mathbin}{AMSb}{"5A}

\begin{document}
\title[Parallel Chip-Firing on the Complete Graph]{Parallel Chip-Firing on the Complete Graph: Devil's Staircase and Poincar\'{e} Rotation Number}
\author[Lionel Levine]{Lionel Levine}
\address{
%\addressmark{1}
Department of Mathematics, Massachusetts Institute of Technology, 77 Massachusetts Ave., Cambridge, MA 02139; {\tt\url{http://math.mit.edu/~levine}}
}
\thanks{The author is supported by a National Science Foundation Postdoctoral Research Fellowship}
\date{January 6, 2010}

\begin{abstract}
We study how parallel chip-firing on the complete graph $K_n$ changes behavior as we vary the total number of chips.  Surprisingly, the activity of the system, defined as the average number of firings
per time step, does not increase smoothly in the number of chips; instead it
remains constant over long intervals, punctuated by sudden jumps.  In
the large $n$ limit we find a ``devil's staircase'' dependence of activity
on the number of chips.  The proof proceeds by reducing the chip-firing
dynamics to iteration of a self-map of the circle $S^1$, in such a way that the
activity of the chip-firing state equals the Poincar\'{e} rotation
number of the circle map.  The stairs of the devil's staircase correspond to periodic chip-firing states of small period.
\end{abstract}

\maketitle

\section{Introduction}

In this paper we explore a connection between the Poincar\'{e} rotation number of a circle map $S^1 \to S^1$ and the behavior of a discrete dynamical system known variously as \emph{parallel chip-firing} \cite{BG,BLS} or the (deterministic) \emph{fixed energy sandpile} \cite{BTW,VDMZ}.  We use this connection to shed light on two intriguing features of parallel chip-firing, \emph{mode locking} and \emph{short period attractors}.  Ever since Bagnoli, Cecconi, Flammini, and Vespignani \cite{BCFV} found evidence of mode locking and short period attractors in numerical experiments in 2003, these two phenomena have called out for a mathematical explanation.
%This paper is intended as a first step toward such an explanation.

In parallel chip-firing on the complete graph $K_n$, each vertex $v \in [n]=\{1,\ldots,n\}$ starts with a pile of $\sigma(v) \geq 0$ chips.  A vertex with $n$ or more chips is \emph{unstable}, and can \emph{fire} by sending one chip to each vertex of $K_n$ (including one chip to itself).  The \emph{parallel update} rule fires all unstable vertices simultaneously, yielding a new chip configuration $U\sigma$ given by
	\begin{equation} \label{parallelupdaterule}
	 U\sigma(v) = \begin{cases} \sigma(v) + r(\sigma), & \sigma(v) \leq n-1 \\
						    \sigma(v) - n + r(\sigma), & \sigma(v) \geq n. 	\end{cases} \end{equation}
Here 
	\[ r(\sigma) = \# \{v| \sigma(v) \geq n\} \]
is the number of unstable vertices.  Write $U^0 \sigma = \sigma$, and $U^t \sigma = U(U^{t-1} \sigma)$ for $t \geq 1$.

Note that the total number of chips in the system is conserved.  In particular, only finitely many different states are reachable from the initial configuration $\sigma$, so the sequence 
%$\sigma, U\sigma, U^2 \sigma, \ldots$
$(U^t \sigma)_{t\geq 0}$ 
is eventually periodic: there exist integers $m \geq 1$ and $t_0 \geq 0$ such that 
	\begin{equation} \label{eventuallyperiodic} U^{t+m} \sigma = U^t \sigma  \qquad \forall t \geq t_0. \end{equation}

The \emph{activity} of $\sigma$ is the limit
		\begin{equation} \label{theactivity} a(\sigma) = \lim_{t \to \infty} \frac{\alpha_t}{nt}. \end{equation}
where
	\[ \alpha_t = \sum_{s=0}^{t-1} r(U^s \sigma) \]
is the total number of firings performed in the first $t$ updates. By (\ref{eventuallyperiodic}), the limit in (\ref{theactivity}) exists and equals
%	$\frac{1}{mn} \sum_{t=t_0}^{t_0+m-1} r(U^t \sigma)$. Since $0 \leq \alpha_t \leq nt$, we have $0 \leq a(\sigma) \leq 1$.
	$\frac{1}{mn} (\alpha_{t+m} - \alpha_{t})$ for any $t \geq t_0$.  
Since $0 \leq \alpha_t \leq nt$, we have $0 \leq a(\sigma) \leq 1$.

\begin{figure}
\begin{center}
\boxed{\includegraphics[width=.46\textwidth]{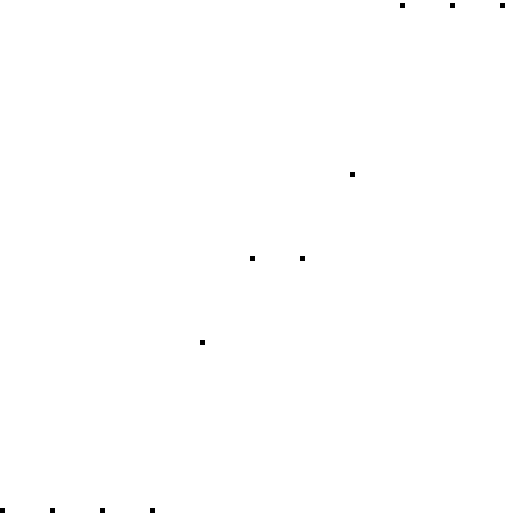}}\hspace{1mm}
\boxed{\includegraphics[width=.46\textwidth]{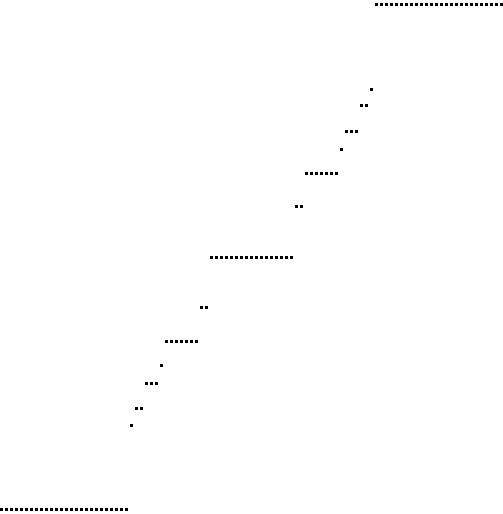}} \\
\vspace{2mm}
\boxed{\includegraphics[width=.46\textwidth]{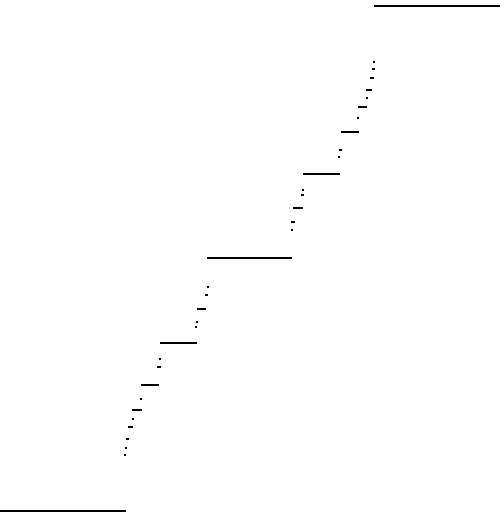}}\hspace{1mm}
\boxed{\includegraphics[width=.46\textwidth]{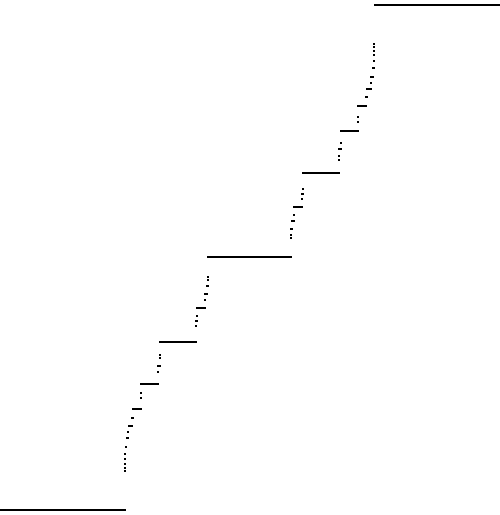}} \\
\end{center}

%\begin{center}
%\begin{tabular}{ccc}
%\includegraphics[scale=.3]{slope2-K_10-fat.png}& ~ &
%\includegraphics[scale=.3]{slope2-K_100-fat.png} \\
%$K_{10}$ & ~ & $K_{100}$
%\end{tabular}
%\end{center}
%~\vspace{4mm}
%\begin{center}
%\begin{tabular}{ccc}
%\includegraphics[scale=.3]{slope2-K_1000-fat.png}& ~ &
%\includegraphics[scale=.3]{slope2-K_10000-fat.png} \\
%$K_{1000}$ & ~ & $K_{10000}$
%\end{tabular}
%\end{center}

\caption{The activity phase diagrams $\tilde{s}_n(k) = a(\sigma_n+k)$, 
for $n=10$ (top left), $100$ (top right), $1000$ (bottom left), and $10000$, where $\sigma_n$ is given by~(\ref{slope2}).  On the  horizontal axis, $k$ runs from $0$ to $n$.  On the vertical axis, $\tilde{s}_n$ runs from $0$ to~$1$.}
\label{slope2staircases}
\end{figure}

Following \cite{BCFV}, we ask: how does the activity change when chips are added to the system?  If $\sigma_n$ is a chip configuration on $K_n$, write $\sigma_n+k$ for the configuration obtained from $\sigma_n$ by adding $k$ chips at each vertex.  The function
	\[ \tilde{s}_n(k) = a(\sigma_n+k) \]
is called the \emph{activity phase diagram} of~$\sigma_n$.  Surprisingly, for many choices of $\sigma_n$, the function $\tilde{s}_n$ looks like a staircase, with long intervals of constancy punctuated by sudden jumps (Figure~\ref{slope2staircases}).  This phenomenon is known as \emph{mode locking}: if the system is in a preferred mode, corresponding to a wide stair in the staircase, then even a relatively large perturbation in the form of adding extra chips will not change the activity.  For a general discussion of mode locking, see \cite{Lagarias}.  

To quantify the idea of mode locking in our setting, suppose we are given an infinite family of chip configurations $\sigma_2, \sigma_3, \ldots$ with $\sigma_n$ defined on $K_n$.  Suppose that $0 \leq \sigma_n(v) \leq n-1$ for all $v \in [n]$, and that for all $0 \leq x \leq 1$
	\begin{equation} \label{converginghistogramsintro} \frac1n \# \{v \in [n] \,|\, \sigma_n(v) < nx\} \to F(x) \qquad \mbox{as } n\to \infty \end{equation}
for a continuous function $F: [0,1] \to [0,1]$.  Then according to Theorem~\ref{mainconvergence}, the activity phase diagrams $\tilde{s}_n$, suitably rescaled, converge pointwise to a continuous, nondecreasing function $s : [0,1] \to [0,1]$.  Moreover, under a mild additional hypothesis, Proposition~\ref{devil'sstaircase} says that this limiting function $s$ is a \emph{devil's staircase}: it is locally constant on an open dense subset of $[0,1]$.  For each rational number $p/q \in [0,1]$ there is a stair of height $p/q$, that is, an interval of positive length on which $s$ is constant and equal to $p/q$.

Related to mode locking, a second feature observed in simulations of parallel chip-firing is non-ergodicity: in trials performed with random initial configurations on the $n \times n$ torus, the activity observed in individual trials differs markedly from the average activity observed over many trials \cite{VDMZ}.  The experiments of \cite{BCFV} suggested a reason: the chip-firing states in locked modes, corresponding to stairs of the devil's staircase, tend to be periodic with very small period. 
%   We study these \emph{short period attractors} in section~\ref{sec:shortperiodattractors}, where we show that all chip configurations on $K_n$ have eventual period at most $n$ (Proposition~\ref{distinctheights}).
 %%, and that the eventual period of $\sigma$ equals the denominator of its activity expressed in lowest terms (Lemma~\ref{denominatorofactivity}).  
    %We also find a small ``window'' in which all configurations have eventual period two (Theorem~\ref{period2window}).   
    We study these \emph{short period attractors} in Theorem~\ref{lotsofsmallperiods}.  Under the same hypotheses that yield a devil's staircase in Propositon~\ref{devil'sstaircase}, for each $q\in \N$, at least a constant fraction $c_q n$ of the states $\{\sigma_n + k\}_{k=0}^n$ have eventual period~$q$.

To illustrate these results, consider the chip configuration  $\sigma_n$ on $K_n$ defined by
	\begin{equation} \label{slope2} \sigma_n(v) = \floor{\frac{n}{4}} + \floor{\frac{v-1}{2}}, \qquad v=1,\ldots,n. \end{equation}
Here $\floor{x}$ denotes the greatest integer $\leq x$.  This family of chip configurations satisfies (\ref{converginghistogramsintro}) with 
	\begin{equation} \label{slope2histogram} F(x) = \begin{cases} 0, & x \leq \frac14 \\
					   2x-\frac12, & \frac14 \leq x \leq \frac34 \\
					   1, & x \geq \frac34. \end{cases} \end{equation}
The activity phase diagrams of $\sigma_n$ for $n=10,100,1000,10000$ are graphed in Figure~\ref{slope2staircases}.  For example, when $n=10$ we have
	\[ (a(\sigma_{10}+k))_{k=0}^{10} = (0,\,0,\,0,\,0,\,1/3,\,1/2,\,1/2,\,2/3,\,1,\,1,\,1) \]
and when $n=100$, we have \vspace{3mm}

\noindent $(a(\sigma_{100}+k))_{k=0}^{100}=$(0, 0, 0, 0, 0, 0, 0, 0, 0, 0, 0, 0, 0, 0, 0, 0, 0, 0, 0, 0, 0, 0, 0, 0, 0, 0, 1/6, 1/5, 1/5, 1/4, 1/4, 1/4, 2/7, 1/3, 1/3, 1/3, 1/3, 1/3, 1/3, 1/3, 2/5, 2/5, 1/2, 1/2, 1/2, 1/2, 1/2, 1/2, 1/2, 1/2, 1/2, 1/2, 1/2, 1/2, 1/2, 1/2, 1/2, 1/2, 1/2, 3/5, 3/5, 2/3, 2/3, 2/3, 2/3, 2/3, 2/3, 2/3, 5/7, 3/4, 3/4, 3/4, 4/5, 4/5, 5/6, 1, 1, 1, 1, 1, 1, 1, 1, 1, 1, 1, 1, 1, 1, 1, 1, 1, 1, 1, 1, 1, 1, 1, 1, 1, 1). \\

As $n$ grows, the denominators of these rational numbers grow remarkably slowly: the largest denominator is~$11$ for $n=1000$, and~$13$ for $n=10000$. 
Moreover, for any fixed $n$ the very smallest denominators occur with greatest frequency. For example, when $n=10000$, there are $1667$ values of $k$ for which $a(\sigma_{n}+k)=\frac12$, and $714$ values of $k$ for which  $a(\sigma_{n}+k)=\frac13$; but for each $p=1,\ldots,12$ there is just one value of $k$ for which $a(\sigma_{n}+k)=\frac{p}{13}$.  
In Lemma~\ref{denominatorofactivity}, we relate these denominators to the periodicity: 
if $a(\sigma)=p/q$ in lowest terms, then $\sigma$ has eventual period $q$.
%In Lemma~\ref{denominatorofactivity}, we show that if $\sigma$ is any chip configuration on $K_n$ that the eventual period of a state equals the denominator of its activity expressed in lowest terms.  
%Thus in the example above, when $n=10000$, all of the states $\{\sigma_n+k\}_{k=0}^n$ are eventually periodic with period at most $13$.  
%This is the phenomenon of \emph{short period attractors} observed in \cite{BCFV}.

The remainder of the paper is organized as follows.  In section~\ref{sec:circlemap} we show how to construct, given a chip configuration~$\sigma$ on~$K_n$, a circle map $f : S^1 \to S^1$ whose rotation number equals the activity of~$\sigma$.  This construction resembles the one-dimensional particle/barrier model of~\cite{JL}.  In section~\ref{sec:staircase} we use the circle map to prove our main results on mode locking, Theorem~\ref{mainconvergence} and Proposition~\ref{devil'sstaircase}.  Short period attractors are studied in section~\ref{sec:shortperiodattractors}, where we show that all states on $K_n$ have eventual period at most $n$ (Proposition~\ref{distinctheights}).  Finally, in Theorem~\ref{period2window}, we find a small ``window'' in which all states have eventual period two.   
%Finally, in Theorem~\ref{lotsofsmallperiods}, under the same hypotheses that yield a devil's staircase in Propositon~\ref{devil'sstaircase}, we show that for each $q\in \N$, a constant fraction $c_q n$ of the states $\{\sigma_n + k\}_{k=0}^n$ have eventual period~$q$.

In parallel chip-firing on a general graph, a vertex is unstable if it has at least as many chips as its degree, and fires by sending one chip to each neighbor; at each time step, all unstable vertices fire simultaneously.  Many questions remain concerning parallel chip-firing on graphs other than $K_n$.  If the underlying graph is a tree \cite{BG} or a cycle \cite{Dall'Asta}, then instead of a devil's staircase of infinitely many preferred modes, there are just three: activity $0$, $\frac12$ and $1$.  On the other hand, the numerical experiments of \cite{BCFV} for parallel chip-firing on the $n\times n$ torus suggest a devil's staircase in the large~$n$ limit.  Our arguments rely quite strongly on the structure of the complete graph, whereas the mode locking phenomenon seems to be widespread.  It would be very interesting to relate parallel chip-firing on other graphs to iteration of a circle map (or perhaps a map on a higher-dimensional manifold) in order to explain the ubiquity of mode locking.
 
\section{Construction of the Circle Map}
\label{sec:circlemap}

We first introduce a framework of generalized chip configurations, which will encompass chip configurations on $K_n$ for all $n$.  To each generalized chip configuration we associate a probability measure on the interval $[0,2)$, and to each such measure we associate a circle map $S^1 \to S^1$.  We will define an update rule $U$ describing the dynamics on each of these objects, so that the following diagram commutes. 
\begin{equation}
\label{thediagram}
\begin{CD}
\mbox{chip configurations on $K_n$} @>U>> \mbox{chip configurations on $K_n$} \\
@VVV @VVV \\
\mbox{generalized chip configurations} @>U>> \mbox{generalized chip configurations} \\
@VVV @VVV \\
\mbox{measures on $[0,2)$} @>U>> \mbox{measures on $[0,2)$} \\
@VVV @VVV \\
\mbox{lifts of circle maps $S^1 \to S^1$} @>U>> \mbox{lifts of circle maps $S^1 \to S^1$}
\end{CD}
\end{equation}
~\\

Write $\Leb$ for Lebesgue measure on $[0,1]$.  A \emph{generalized chip configuration} is a measurable function
	\[ \eta : [0,1] \to [0,\infty). \]
Let
	\[ r(\eta) = \Leb \{ x | \eta(x) \geq 1 \} \]	
and define the update rule
	\begin{equation} \label{generalizedparallelupdate} 
		U\eta (x) = \begin{cases} 
			\eta(x) + r(\eta), &\eta(x) < 1 \\
			\eta(x) - 1 + r(\eta), &\eta(x) \geq 1. 
		\end{cases}
	\end{equation}
If $\sigma$ is a chip configuration on $K_n$, we define its associated generalized chip configuration $\psi(\sigma)$ by
%	\[ \eta(\sigma) = \frac{1}{n} \sum_{v=1}^n \sigma(v) 1_{[\frac{v-1}{n},\frac{v}{n})} \]
%% constant on intervals: not good because the resulting circle map is discontinuous
	\begin{equation} \label{piecewiselinear} \psi(\sigma)(x) = \frac{\sigma( \ceiling{nx} ) + \ceiling{nx}-nx}{n}   \end{equation}
% linear on intervals is better.
where $\ceiling{y}$ denotes the least integer $\geq y$.  Our first lemma checks that the top square of (\ref{thediagram}) commutes.

\begin{lemma} 
\label{intertwining}
$U \circ \psi = \psi \circ U$.
\end{lemma} 

\begin{proof}
Let $\sigma$ be a chip configuration on $K_n$, and write $\eta = \psi(\sigma)$.
For $v\in [n]$, if $\sigma(v) \geq n$, then $\eta(x) \geq 1$ for all $x \in (\frac{v-1}{n},\frac{v}{n}]$.  Conversely, if $\sigma(v) \leq n-1$, then $\eta(x) <1$ for all $x \in (\frac{v-1}{n},\frac{v}{n}]$.  Hence $r(\eta) = r(\sigma)/n$.  

Fix $x \in [0,1]$ and let $v=\ceiling{nx}$.  Then
	\begin{align*} U(\psi(\sigma))(x) &= \eta(x) - 1_{\{\eta(x)\geq 1\}} + r(\eta) \\
					&= \frac{\sigma(v)+v-nx}{n} - 1_{\{\sigma(v)\geq n\}} + \frac{r(\sigma)}{n} \\
					&= \frac{U\sigma(v)+v-nx}{n} \\
					&= \psi(U\sigma)(x). \qed \end{align*}
\renewcommand{\qedsymbol}{}
\end{proof}
 
We remark that Lemma~\ref{intertwining} would hold also if we defined $\psi$ using the piecewise constant interpolation $\frac{\sigma(\ceiling{nx})}{n}$.  We have chosen the piecewise linear interpolation~(\ref{piecewiselinear}) because it is better suited to our construction of the circle map, below: the circle map associated to the piecewise linear interpolation will always be continuous.
	
Next we observe a simple consequence of the update rule (\ref{generalizedparallelupdate}), which will allow us to focus on a subset of generalized chip configurations which we call ``confined.''
 
\begin{lemma}
\label{confinedrange}
Let $\eta$ be a generalized chip configuration, and let $x \in [0,1]$.  If~$\eta(x) < 2$, then 
	\[ r(\eta) \leq U\eta(x) < 1+r(\eta). \]
\end{lemma}

\begin{proof}
The first inequality is immediate from the definition of $U\eta$.
For the second inequality, if $\eta(x)<1$, then 
	\[ U\eta(x) = \eta(x) + r(\eta) < 1+ r(\eta) \]
while if $1\leq \eta(x)<2$, then
	\[ U\eta(x) = \eta(x) - 1 + r(\eta) < 1 + r(\eta). \qed \]
\renewcommand{\qedsymbol}{}
\end{proof}

\begin{definition}
A generalized chip configuration $\eta$ is \emph{preconfined} if it satisfies
	\begin{enumerate}
	\item[(i)] $\eta(x) <2$ for all $x\in [0,1]$.
	\end{enumerate}
If, in addition, 
there exists $r \in [0,1]$ such that
	\begin{enumerate}
	\item[(ii)] $r \leq \eta(x) < 1+r$ for all $x \in [0,1]$
	\end{enumerate}
then $\eta$ is \emph{confined}.
\end{definition}

By Lemma~\ref{confinedrange}, if $\eta$ is preconfined, then $U\eta$ is confined.

Note that from (\ref{parallelupdaterule})
	\[ U\eta(x) \equiv \eta(x)+r(\eta)  \modulo{1}. \]
Iterating yields the congruence
	\begin{equation} \label{iteratedcongruence}  U^t \eta (x) \equiv \eta(x) + \beta_t \modulo{1} \end{equation}
where
	\[ \beta_t = \sum_{s=0}^{t-1} r(U^s \eta). \]
	
Next we find a recurrence for the sequence $\beta_t$.  We start with the following lemma.
%Our next lemma characterizes the vertices that fire at time $t+1$.

\begin{lemma}
\label{whichsitesareunstable}
%If $\sigma$ is confined, 
If $U^t \eta$ is preconfined, then $U^{t+1}\eta(x) \geq 1$ if and only if 
	\[ \eta(x) \in \Z - (\beta_t, \beta_{t+1}]. \]
%	 \[ \sigma(x) \equiv -y~(\mbox{\em mod } 1) \] 
%for some $\alpha_t < y \leq \alpha_{t+1}$.
\end{lemma}

\begin{proof}
%By Lemma~\ref{transientheights} we have $U^t \sigma \leq 2n-1$, so 
By Lemma~\ref{confinedrange}, since $U^t \eta$ is preconfined, we have	
	\[ r_t \leq U^{t+1} \eta(x) < 1+r_t \]
for all $x \in [0,1]$, where $r_t = r(U^t \eta)$.
Thus $U^{t+1} \eta(x) \geq 1$ if and only if 
	\[ U^{t+1} \eta (x) \in \Z + [0,r_t).  \] 
By (\ref{iteratedcongruence}), this condition is equivalent to
	\[ \eta(x) + \beta_{t+1} \in \Z + [0,r_t). \]
Using the fact that $\beta_{t+1} = \beta_t + r_t$, this in turn is equivalent to
	\[ \eta(x) \in \Z + [-\beta_{t+1}, r_t - \beta_{t+1}) = \Z - (\beta_t, \beta_{t+1}]. \qed \]
\renewcommand{\qedsymbol}{}	
\end{proof}

The essential information in a generalized chip configuration $\eta$ is contained in the associated probability measure $\mu$ on $[0,\infty)$ given by
	\[ \mu(A) = \Leb (\eta^{-1}(A)) \]
for Borel sets $A \subset [0,\infty)$.	
If $\eta$ arises from a chip configuration $\sigma$ on $K_n$, then $\mu([a,b])$ is the proportion of vertices that have between $an$ and $bn$ chips.
For~$y\in \R$, write $\xi^y \mu$ for the translated measure $\xi^y \mu (A) = \mu(A-y)$.  For~$Y \subset \R$, write $\mu|_Y$ for the restricted measure $\mu|_Y(A) = \mu(A\cap Y)$.   Then the update rule (\ref{generalizedparallelupdate}) takes the form
	\begin{equation} \label{measureparallelupdate} U\mu = \xi^{\mu[1,\infty)} (\mu|_{[0,1)}) + \xi^{\mu[1,\infty) -1} (\mu|_{[1,\infty)}). \end{equation}
%	\[ U\mu = \xi^{\mu[1,\infty)} \left( \mu|_{[0,1)} + \xi^{-1} (\mu|_{[1,\infty)}) \right). \]
%Whether one works with generalized chip configurations $\eta$ or their associated measures $\mu$ is largely a matter of convenience; we will usually choose to work with generalized chip configurations.

Now consider the measure $\nu$ on $\R$ given by
	\begin{equation} \label{thisisnu} \nu(A) = \sum_{m \in \Z} \mu(m-A) = \sum_{m \in \Z} \Leb \{x \,|\, m-\eta(x) \in A \}. \end{equation}
%	\begin{equation} \label{thisisnu} \nu(A) = \mu(\Z-A) = 
%		\Leb \{x \,|\, \eta(x) \in \Z-A\}. \end{equation}
By Lemma~\ref{whichsitesareunstable}, if $U^t \eta$ is preconfined, then 
	\[r(U^{t+1}\eta) = \mu(\Z - (\beta_t,\beta_{t+1}]) =  \nu(\beta_t,\beta_{t+1}], \] 
hence
	\begin{equation} \label{thethreetermrecurrence} \beta_{t+2} = \beta_{t+1} + \nu(\beta_t,\beta_{t+1}]. \end{equation}
This gives a recurrence relating three consecutive terms of the sequence~$\beta_t$.  Our next lemma simplifies this recurrence to one relating just two consecutive terms.

\begin{lemma}
\label{twotermrecurrence}
If $\eta$ is preconfined, then for all $t \geq 0$
	\[ \beta_{t+1}=f(\beta_t),\] 
where $f : [0,\infty) \to [0,\infty)$ is given by
	\begin{equation} \label{precursortothecirclemap} f(y) = \beta_1 + \nu(0,y]
	\end{equation}
and $\nu$ is given by (\ref{thisisnu}). 
\end{lemma}

\begin{proof} 
By Lemma~\ref{confinedrange}, $U^t \eta$ is preconfined for all $t\geq 0$, so
from (\ref{thethreetermrecurrence}),
	\begin{align*} 
	\beta_{t+1} - \beta_1 &= \sum_{s=0}^{t-1}  (\beta_{s+2}-\beta_{s+1})  \\
		&=  \sum_{s=0}^{t-1} \nu(\beta_{s},\beta_{s+1}] \\
		&= \nu(0,\beta_t].
	\end{align*}
Hence $\beta_{t+1} = f(\beta_t)$.
\end{proof}

The function $f$ appearing in Lemma~\ref{twotermrecurrence} satisfies
	\begin{align*} f(y+1) &= f(y) + \nu(y,y+1] \\ &= f(y) + 1 \end{align*}
for all $y \geq 0$.  Thus it has a unique extension to a function $f : \R \to \R$ satisfying $f(y+1)=f(y)+1$ for all $y \in \R$.  

Note that $f$ is nondecreasing.  If it is also continuous, 
% is continuity really needed?
then it has a well-defined \emph{Poincar\'e rotation number} \cite{HK,1dim}
	\[ \rot(f) = \lim_{t\to\infty} \frac{f^t(y)}{t} \]
which does not depend on $y$.  Here $f^t$ denotes the $t$-fold iterate $f^t(y) = f(f^{t-1}(y))$, with $f^0=Id$.  

Viewing the circle $S^1$ as $\R/\Z$, the map $f$ descends to a circle map $ \bar{f} : S^1 \to S^1$.  The rotation number $\rot(f)$ measures the asymptotic rate at which the sequence \[ y,\, \bar{f}(y), \, \bar{f}(\bar{f}(y)), \, \ldots \] winds around the circle.

We define the \emph{activity} of a generalized chip configuration $\eta$ as the limit
	\[ a(\eta) = \lim_{t \to \infty} \frac{\beta_t}{t}. \]
From Lemma~\ref{twotermrecurrence} we see that if $\eta$ is preconfined and $f$ is continuous, then this limit exists and equals the rotation number of $f$.

\begin{lemma}
\label{theiterates}
If $\eta$ is preconfined, then $\beta_t = f^t(0)$ for all $t \geq 0$.
\end{lemma}

\begin{lemma}
\label{activityequalsrotationnumber}
If $\eta$ is preconfined and $f$ is continuous, then $a(\eta)=\rho(f)$.
\end{lemma}

The conditions that $\eta$ is preconfined and $f$ is continuous are most succinctly expressed in terms of the associated probability measure:  $\mu$ is supported on $[0,2)$ and is \emph{non-atomic}, that is, $\mu(\{y\}) =0$ for all $y \in [0,2)$.

%Note that the map $f$ is defined in terms of $\beta_1$ and $\nu$, both of which are easily read off from $\eta$.  So given a preconfined configuration~$\eta$, equation~(\ref{precursortothecirclemap}) gives an explicit recipe for writing down a circle map $f$ whose rotation number is the activity of $\eta$.

To complete the commutative diagram (\ref{thediagram}), it remains to describe the dynamics $U$ on lifts of circle maps.  Call a map $f: \R \to \R$ a \emph{monotone degree one lift} if $f$ is continuous, nondecreasing and satisfies
	\begin{equation} \label{degreeonelift} f(y+1) = f(y) + 1 \end{equation}
for all $y \in \R$.  Let $\mathcal{F}$ be the set of monotone degree-one lifts 
$f : \R \to \R$, 
equipped with the $L^\infty$ topology.  We define the update rule $U: \mathcal{F} \to \mathcal{F}$ by
	\[ Uf = R_{-f(0)} \circ f \circ R_{f(0)} \]
%		\[ Uf(x) = f(x+f(0))-f(0). \]
where for $y \in \R$, we write $R_y$ for the translation $x \mapsto x+y$.

%A measure $\mu$ on $[0,\infty)$ is \emph{non-atomic} if $\mu(\{x\})=0$ for all $x \geq 0$.
Let~$\mathcal{M}$ be the space of non-atomic probability measures $\mu$ on the interval~$[0,2)$, equipped with the topology of weak-$*$ convergence: $\mu_n \to \mu$ if and only if $\mu_n[0,y) \to \mu[0,y)$ for all $y\in [0,2)$.  For $\mu \in \mathcal{M}$, define $\phi(\mu) : \R \to \R$ by
	\[ \phi(\mu)(y) = \mu[1,2) + \sum_{m \in \Z} \mu(m-y,m] \]
	%	= \mu[1,2) + \mu(\Z-[0,y)). \]
where if $b<a$ we define $\mu(a,b]:=-\mu[b,a)$.
	
\begin{theorem}
\label{thebottomsquare}
If $\mu \in \mathcal{M}$, then
$\phi(\mu) \in \mathcal{F}$.  The map
$\phi: \mathcal{M} \to \mathcal{F}$ is continuous and preserves the dynamics: $U \circ \phi = \phi \circ U$.
\end{theorem}
\begin{proof}
$\phi(\mu)$ is clearly nondecreasing in $y$, and since $\mu$ is non-atomic, $\phi(\mu)$ is continuous.  Moreover, for any $y \in \R$,
	\[ \phi(\mu)(y+1) - \phi(\mu)(y) = \sum_{m \in \Z} \mu(m-y-1,m-y] = \mu(\R)=1 \]
so $\phi(\mu) \in \mathcal{F}$.

If $\mu_n \in \mathcal{M}$ and $\mu_n \to \mu \in \mathcal{M}$, then the cumulative distribution functions
	$F_n(y) = \mu_n[0,y)$
converge pointwise to $F(y):=\mu[0,y)$.  Since $F_n$ and $F$ are continuous and nondecreasing, this convergence is uniform in~$y$.  Since $\mu_n$ and $\mu$ are supported on the interval $[0,2)$, we have for $y\in [0,1]$
	\begin{align*} \phi(\mu_n)(y) &= \mu_n[1,2) + \mu_n(1-y,1] + \mu_n(2-y,2] \\
						&= 2F_n(2) - F_n(1-y) - F_n(2-y). \end{align*}
The right side converges uniformly in $y$ to $2F(2)-F(1-y)-F(2-y)=\phi(\mu)(y)$.  Since $\phi(\mu_n)$ and $\phi(\mu)$ are degree one lifts, we have
	\[ \sup_{y \in \R} | \phi(\mu_n)(y) - \phi(\mu)(y) | = \sup_{y \in [0,1]} | \phi(\mu_n)(y) - \phi(\mu)(y) | \to 0 \]
as $n \to \infty$, which shows that $\phi$ is continuous.

Let $\beta = \mu[1,2) = \phi(\mu)(0)$.  Then
	\begin{align*} (U\circ \phi)(\mu)(y) &= (R_{-\beta} \circ \phi(\mu) \circ R_{\beta})(y) \\
					&= \sum_{m \in \Z} \mu(m-y-\beta,m]. \end{align*}
Write $\mu=\mu_0 +\mu_1$, where $\mu_0 = \mu|_{[0,1)}$ and $\mu_1 = \mu|_{[1,2)}$.  By (\ref{measureparallelupdate}), we have
	\[ U\mu = \xi^\beta \mu_0 + \xi^{\beta-1} \mu_1 \]
hence
	\begin{align*} U\mu[1,2) &= \mu_0 [1-\beta,2-\beta) + \mu_1 [2-\beta,3-\beta) \\					&= \mu[1-\beta,1) + \mu[2-\beta,2) \\
					&= \sum_{m \in \Z} \mu(m-\beta,m] \end{align*}
where in the last line we have used that $\mu$ is non-atomic and supported on $[0,2)$.  Moreover	
	\begin{align*}
	 \sum_{m \in \Z} U\mu(m-y,m] &= 
		\sum_{m \in \Z} \mu_0(m-y-\beta,m-\beta] \,+ \\ 
		&\qquad\qquad	+\sum_{m\in \Z} \mu_1(m-y+1-\beta,m+1-\beta]  \\
					&= \sum_{m \in \Z} \mu(m-y-\beta,m-\beta].
	\end{align*}
We conclude that
	\begin{align*} (\phi\circ U)(\mu)(y) &= U\mu[1,2) + \sum_{m \in \Z} U\mu(m-y,m] \\
					&= \sum_{m \in \Z} \mu(m-y-\beta,m] \end{align*}
so $\phi \circ U = U \circ \phi$.
\end{proof}

One naturally wonders how to generalize the construction described in this section to chip-firing on graphs other than $K_n$.  A key step may involve identifying invariants of the dynamics.  On $K_n$, these invariants take a very simple form: by (\ref{iteratedcongruence}), for any two vertices $v,w \in [n]$, the difference
	\[ U^t\sigma(v) - U^t\sigma(w) \quad \mod n \]
does not depend on $t$.  Analogous invariants for parallel chip-firing on the $n\times n$ torus are classified in \cite{CDVV}, following the approach of \cite{Dhar}.

\section{Devil's Staircase}
\label{sec:staircase}

  Let $f,f_n,g$ be monotone degree one lifts (\ref{degreeonelift}), and denote by $\bar{f}, \bar{f}_n, \bar{g}$ the corresponding circle maps $S^1 \to S^1$.  Write $f \leq g$ if $f(x) \leq g(x)$ for all $x \in \R$, and $f<g$ if $f(x) < g(x)$ for all $x \in \R$.  We will make use of the following well-known properties of the rotation number.  For their proofs, see, for example \cite{HK,1dim}.

\begin{itemize}
\item \textbf{Monotonicity}. If $f \leq g$, then $\rho(f) \leq \rho(g)$.

\item \textbf{Continuity}. If $\sup |f_n - f| \to 0$, then $\rho(f_n) \to \rho(f)$.

\item \textbf{Conjugation Invariance}. If $g$ is a homeomorphism, then $\rho(g \circ f \circ g^{-1}) = \rho(f)$.
% note: this is only true for increasing homeomorphisms, but degree one implies increasing.

\item \textbf{Instability of an irrational rotation number}. If $\rho(f) \notin \Q$, and $f_1 < f < f_2$, then $\rho(f_1) < \rho(f) < \rho(f_2)$.

\item \textbf{Stability of a rational rotation number}. If $\rho(f) =p/q \in \Q$, and $\bar{f}^q \neq Id : S^1 \to S^1$, then for sufficiently small $\epsilon>0$, either 
	\[ \rho(g) = p/q \mbox{ whenever } f \leq g \leq f+\epsilon, \] 
or 
	\[ \rho(g)=p/q  \mbox{ whenever } f-\epsilon \leq g \leq f. \]

\end{itemize}

Let $\sigma_2, \sigma_3, \ldots$ be a sequence of chip configurations, with $\sigma_n$ 
defined on~$K_n$.  Suppose $\sigma_n$ is stable, i.e., 
	\begin{equation} \label{stability} 0 \leq \sigma_n(v) \leq n-1 \end{equation}
for all $v \in [n]$.  Moreover, suppose that there is a continuous function $F: [0,1] \to [0,1]$, such that for all $0 \leq x \leq 1$
	\begin{equation} \label{converginghistograms} \frac1n \# \{v \in [n] \,|\, \sigma_n(v) < nx\} \to F(x) \end{equation}
% can we let F be discontinuous and require this to hold only at continuity points x?
%uniformly in $x$,  (automatic since F is continuous and F_n is nondecreasing)
as $n \to \infty$.  Since the left side is nondecreasing in $x$, this convergence is uniform in $x$.  We define $\Phi: \R \to \R$ by
	\begin{equation} \label{extensiontoR} \Phi(x) = \ceiling{x} - F(\ceiling{x}-x). \end{equation}
%	\begin{equation} \label{extensiontoR} \tilde F(x)= \floor{x} + F(x-\floor{x}). \end{equation}
%	\begin{equation} \label{extensiontoR} \tilde F(x+m)= \tilde F(x)+m, \qquad m \in \Z. \end{equation}
Note that (\ref{stability}) and (\ref{converginghistograms}) force $F(0)=0$ and $F(1)=1$; by compactness, $F$ is uniformly continuous on $[0,1]$, and hence $\Phi$
%its extension (\ref{extensiontoR}) 
is uniformly continuous on $\R$.

\begin{figure}
\begin{center}
\begin{tabular}{ccc}
\includegraphics[width=.305\textwidth]{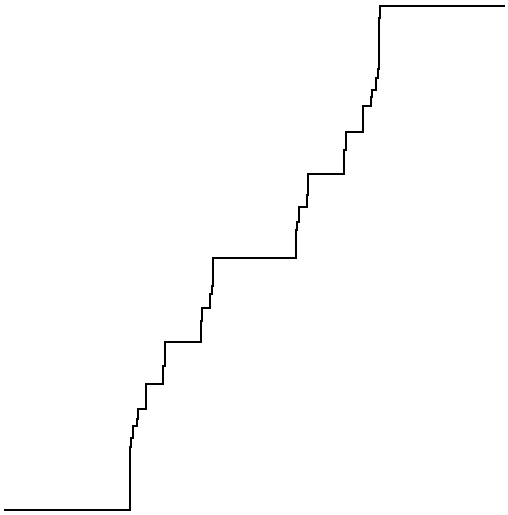} &
\includegraphics[width=.305\textwidth]{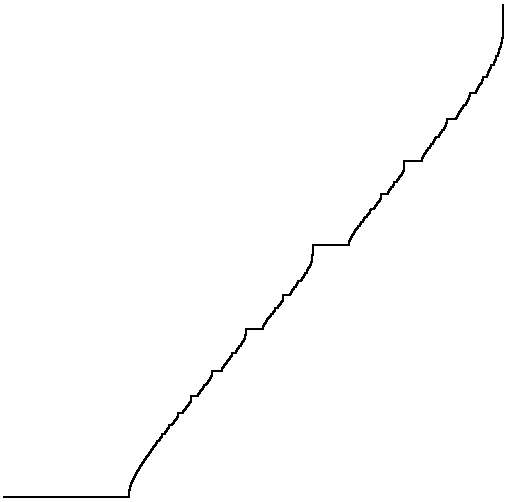} &
\includegraphics[width=.305\textwidth]{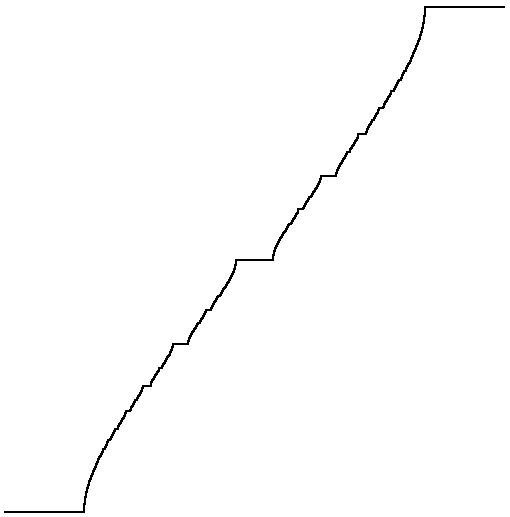} \\
(a) & (b) & (c)
\end{tabular}
\end{center}
\caption{The devil's staircase $s(y)$, when (a) $F(x)$ is given by~(\ref{slope2histogram}); (b) $F(x)=\sqrt{x}$ for $x\in [0,1]$; and (c) $F(x) = x + \frac{1}{2\pi} \sin 2\pi x$.  On the horizontal axis $y$ runs from $0$ to $1$, and on the vertical axis $s(y)$ runs from $0$ to~$1$.
 }
\label{samplestaircases}
\end{figure}

The \emph{rescaled activity phase diagram} of $\sigma_n$ is the function $s_n : [0,1] \to [0,1]$ defined by
	\[ s_n(y) = a(\sigma_n + ny). \]
As $n\to \infty$, the $s_n$ have a pointwise limit identified in our next result. 

\begin{theorem}
\label{mainconvergence}
If (\ref{stability}) and (\ref{converginghistograms}) hold, then for each $y \in [0,1]$ we have
 	\[ s_n(y) \to s(y) := \rot(R_y \circ \Phi) \]
as $n \to \infty$, where $\Phi$ is given by (\ref{extensiontoR}), and $R_y(x) = x+y$.
\end{theorem}

In Proposition~\ref{devil'sstaircase}, below, we show that under an additional mild hypothesis, the limiting function $s(y)$ is a devil's staircase.  Examples of these staircases for different choices of $F$ are shown in Figure~\ref{samplestaircases}.

To prepare for the proof of Theorem~\ref{mainconvergence}, let $\eta_{n,y}$ be the generalized chip configuration associated to $\sigma_n + ny$, let $\mu_{n,y}$ be the associated measure on~$[0,2)$, and let $f_{n,y} = \phi(\mu_{n,y})$ be the corresponding circle map lift.

\begin{lemma}
\label{theuglystuff}
Let $y \in [0,1]$.
If (\ref{stability}) and (\ref{converginghistograms}) hold, then as $n \to \infty$
	\[ \sup |f_{n,y} - \Phi \circ R_y| \to 0 \]
uniformly in $y$.
%	\[ f_{n,y}(x) \to \Phi(x+y) \]
% uniformly in x and y. 
\end{lemma}

\begin{proof}
The measures $\mu_{n,0}$ are nonatomic, and by (\ref{stability}) they are supported on~$[0,1)$, so for $x \in [0,1]$ we have
	\begin{align*} f_{n,y}(x-y) 
	%&= \phi(\xi^y \mu_{n,0})(x-y) \\
		&= \mu_{n,y}[1,2) + \sum_{m \in \Z} \mu_{n,y}(m-x+y,m] \\
		&= \mu_{n,0}[1-y,2-y) + \sum_{m \in \Z} \mu_{n,0}(m-x,m-y] \\
		&= \mu_{n,0}[1-y,1) + \mu_{n,0}(1-x,1-y] \\
		&= \mu_{n,0}[1-x,1) \\
		&= f_{n,0}(x).
	\end{align*}
Since $f_{n,y}$ and $f_{n,0} \circ R_y$ are degree one lifts that agree on $[-y,1-y]$, they agree everywhere.   

By (\ref{converginghistograms}) we have $\mu_{n,0} \to \mu$, where $\mu([a,b]) = F(b)-F(a)$.  Since $\mu$ is supported on $[0,1)$, we have for $x \in [0,1]$
	\[ \phi(\mu)(x) = \mu(1-x,1] = 1 - F(1-x) = \Phi(x). \]
Since $\Phi$ and $\phi(\mu)$ are degree one lifts that agree on $[0,1]$, they agree everywhere.  Hence
	\begin{align*} \sup |f_{n,y} - \Phi \circ R_y | 
	&= \sup | (f_{n,0} - \Phi) \circ R_y | \\
	&= \sup |f_{n,0} - \Phi |  \\
	&= \sup |\phi(\mu_{n,0}) - \phi(\mu)|.
	\end{align*}
The right side does not depend on $y$, and tends to zero as $n \to \infty$ by the continuity of $\phi$ (Theorem~\ref{thebottomsquare}).
\end{proof}

\begin{proof}[Proof of Theorem~\ref{mainconvergence}]
By Lemma~\ref{activityequalsrotationnumber}, Lemma~\ref{theuglystuff} and the continuity of the rotation number,
	\[ s_n(y) = a(\eta_{n,y}) = \rho(f_{n,y}) \to \rho (\Phi \circ R_y) \]
By the conjugation invariance of the rotation number, $\rho(\Phi \circ R_y) = \rho(R_y \circ \Phi)$, which completes the proof.
\end{proof}

Write $\Phi_y = R_y \circ \Phi$, and let $\bar{\Phi}_y : S^1 \to S^1$ be the corresponding circle map.  We will call a function $s : [0,1] \to [0,1]$ a \emph{devil's staircase} if it is continuous, nondecreasing, nonconstant, and locally constant on an open dense set.  Next we show that if
\begin{equation}
\label{notidentity}
 \left(\bar{\Phi}_y\right)^q \neq Id \quad \mbox{for all $y\in [0,1)$ and all $q \in \N$,}
\end{equation}
then the limiting function $s(y)$ in Theorem~\ref{mainconvergence} is a devil's staircase.

\begin{prop}
\label{devil'sstaircase}
The function $s(y)=\rho(\Phi_y)$ continuous and nondecreasing in $y$.  If $z \in [0,1]$ is irrational, then $s^{-1}(z)$ is a point.  Moreover, if 
%$(\Phi_y)^m \neq Id$ for all $y\in S^1$ and all positive integers $m$, 
(\ref{notidentity}) holds, then for every rational number $p/q \in [0,1]$ the fiber $s^{-1}(p/q)$ is an interval of positive length.
\end{prop}

Substantially similar results appear in \cite{HK,1dim}; we include a proof here for the sake of completeness.

\begin{proof}
The monotonicity of the rotation number implies that $s$ is nondecreasing.  Since $\sup |\Phi_y - \Phi_{y'} | = |y-y'|$,
the continuity of the rotation number implies $s$ is continuous.

Since $\Phi(m)=m$ for all $m\in \Z$, we have $s(0)=0$ and $s(1)=1$.  By the intermediate value theorem, $s: [0,1] \to [0,1]$ is onto.  If $z=s(y)$ is irrational, and $y_1 < y < y_2$, then $\Phi_{y_1} < \Phi_y < \Phi_{y_2}$, hence
	\[ s(y_1) < s(y) < s(y_2) \]
by the instability of an irrational rotation number.  It follows that $s^{-1}(z) = \{y\}$.

If $s(y)=p/q \in \Q$, then by the stability of a rational rotation number, since $(\bar{\Phi}_y)^q \neq Id$, there exists an interval $I$ of positive length (take either $I = [y-\epsilon,y]$ or $I=[y,y+\epsilon]$ for small enough $\epsilon$) such that $\rho(\Phi_{y'}) = p/q$ for all $y' \in I$.  Hence $I \subset s^{-1}(p/q)$.
\end{proof}

%We remark that condition (\ref{notidentity}) is a rather weak restriction.  For example, it holds whenever $\bar{\Phi}$ is a Morse-Smale diffeomorphism \cite{1dim}.

Note that if $y \leq y'$, then $f_{n,y} \leq f_{n,y'}$, so $s_n(y) = \rho(f_{n,y})$ is nondecreasing in $y$.  By the continuity of $s$, it follows that the convergence in Theorem~\ref{mainconvergence} is uniform in $y$.  

Our next result shows that in the interiors of the stairs, we in fact have $s_n(y)=s(y)$ for sufficiently large $n$.

\begin{prop}
\label{strongerconvergence}
Suppose that (\ref{stability}), (\ref{converginghistograms}) and (\ref{notidentity}) hold.  If $s^{-1}(p/q) = [a,b]$, then for any $\epsilon>0$
	\[ [a+\epsilon, b-\epsilon] \subset s_n^{-1}(p/q) \]
for all sufficiently large $n$.
% formally weaker, but equivalent by monotonicity: 
% if $y$ is in the interior of the interval $s^{-1}(p/q)$, then 
%	\[ s_n(y) = p/q \]
%for all sufficiently large $n$.
\end{prop}

\begin{proof}
%Since $\Phi$ is uniformly continuous on $[0,1]$, there exists $\delta>0$ such that
%	\[ \sup | \Phi \circ R_y - \Phi \circ R_{y'} | < \epsilon \]
%whenever $|y-y'|<\delta$.  
%not needed I think
Fix $\epsilon>0$.  By Lemma~\ref{theuglystuff}, for sufficiently large $n$ we have 
	\[ \sup |f_{n,y} - \Phi \circ R_y| < \epsilon \]
for all $y \in [0,1]$.  Then $f_{n,y} \geq R_{-\epsilon} \circ \Phi \circ R_{y}$,
%$f_{n,a+\epsilon} \geq \Phi \circ R_{a+\epsilon} - \epsilon$ and
%	\[ R_{a+\epsilon} \circ (\Phi \circ R_{a+\epsilon} - \epsilon) \circ R_{a+\epsilon}^{-1} = R_a \circ \Phi \]
so for any $y \geq a+\epsilon$ we have by the monotonicity and conjugation invariance of the rotation number
	\begin{align*} s_n (y) = \rot (f_{n,y}) &\geq \rot (R_{-\epsilon} \circ \Phi \circ R_{y}) \\ 
	&= \rot (R_{y-\epsilon} \circ \Phi) = s(y-\epsilon) \geq s(a) = p/q. \end{align*}

Likewise $f_{n,y} \leq R_{\epsilon} \circ \Phi \circ R_{y}$, hence for any $y \leq b-\epsilon$
	\begin{align*} s_n (y) = \rot (f_{n,y}) &\leq \rot (R_{\epsilon} \circ \Phi \circ R_{y}) \\
	&= \rot (R_{y+\epsilon} \circ \Phi) = s(y+\epsilon) \leq s(b) = p/q. \qed \end{align*}
\renewcommand{\qedsymbol}{}
\end{proof}

\section{Short Period Attractors}
\label{sec:shortperiodattractors}

In this section we explore the prevalance of parallel chip-firing states on~$K_n$ with small period.  

\begin{definition}
A chip configuration $\sigma$ on $K_n$ is \emph{preconfined} if it satisfies
	\begin{enumerate}
	\item[(i)] $\sigma(v) \leq 2n-1$ for all vertices $v$ of $K_n$.
	\end{enumerate}
If, in addition, $\sigma$ satisfies
	\begin{enumerate}
	\item[(ii)] $\max_v \sigma(v) - \min_v \sigma(v) \leq n-1$
	\end{enumerate}
then $\sigma$ is \emph{confined}.
\end{definition}

Equivalently, $\sigma$ is preconfined (confined) if and only if the generalized chip configuration $\eta(\sigma)$ is preconfined (confined) as defined in section~\ref{sec:circlemap}.

Recall (\ref{eventuallyperiodic}) that for any chip configuration~$\sigma$ on~$K_n$, the sequence $(U^t \sigma)_{t\geq 0}$ is eventually periodic.  Denote its transient length by $t_0$; that is
	\[ t_0 = \min \{t\geq 0 \,|\, U^t\sigma = U^{t'} \sigma \mbox{ for some } t'>t\}. \]  
We say that ``$v$ fires at time $t$'' if $U^t \sigma(v) \geq n$.

\begin{lemma}
\label{eventuallyconfined}
If $a(\sigma)<1$, then $U^t \sigma$ is confined for all $t \geq t_0$.
%If $a(\sigma)<1$, then $U^t \sigma (v)\leq 2n-2$ for all $v$ and all $t \geq t_0$.
\end{lemma}

\begin{proof}
If $a(\sigma)<1$, then at each time step, some vertex does not fire.   If a given vertex $v$ fires at time $t$, then
	\[ U^{t+1} \sigma(v) < U^{t} \sigma(v). \]
Since $U^t \sigma(v) \in \Z_{\geq 0}$ for all $t$, there is some time time $T_v$ at which $v$ does not fire.  By Lemma~\ref{confinedrange}, we have $U^t \sigma(v) \leq 2n-1$ for all $t \geq T_v$, and $U^t\sigma$ is confined for all $t>\max_v T_v$.  

For any $t \geq t_0$ we have $U^t \sigma = U^{t'} \sigma$ for infinitely many values of $t'$, hence~$U^{t} \sigma$ is confined.
\end{proof}

For a chip configuration $\sigma$ on $K_n$ and a vertex $v\in [n]$, let
	\[ u_t(\sigma, v) = \# \{0\leq s<t \,|\, U^s \sigma(v) \geq n \} \]
be the number of times $v$ fires during the first $t$ updates.  During these updates, the vertex $v$ emits a total of $n\, u_t(\sigma,v)$ chips and receives a total of $\alpha_t = \sum_w u_t(\sigma,w)$ chips, so that
	\begin{equation} \label{laplacianofodometer} U^t \sigma(v) - \sigma(v) = \alpha_t - n\, u_t(\sigma,v). \end{equation}
An easy consequence that we will use repeatedly is the following.

\begin{lemma}
\label{constantodometer}
A chip configuration $\sigma$ on $K_n$ satisfies
$U^t \sigma = \sigma$ if and only if 
	\begin{equation} \label{eq:constantodometer} u_t(\sigma,v) = u_t(\sigma,w) \end{equation}
for all vertices $v$ and $w$.
\end{lemma}

\begin{proof}
If (\ref{eq:constantodometer}) holds, then $u_t(\sigma,v) = \alpha_t/n$ for all $v$, so $U^t \sigma(v) = \sigma(v)$ by (\ref{laplacianofodometer}).
%If $u_t(\sigma,v) = \ell$ for all vertices $v$, then during the first $t$ updates, each vertex $v$ emits a total of $n\ell$ chips and receives a total of $n\ell$ chips, so $U^t \sigma(v) = \sigma(v)$.
Conversely, if $U^t \sigma = \sigma$, then the left side of (\ref{laplacianofodometer}) vanishes, so $u_t(\sigma,v) = \alpha_t/n$ for all vertices $v$.
\end{proof}

According to our next lemma, if $\sigma$ is confined, then $u_t(\sigma,v)$ and $u_t(\sigma,w)$ differ by at most one.

\begin{lemma}
\label{interlacing}
If $\sigma$ is confined, and $\sigma(v) \leq \sigma(w)$, then for all $t\geq 0$
	\[ u_t(\sigma,v) \leq u_t(\sigma,w) \leq u_t(\sigma,v)+1. \]
\end{lemma}

\begin{proof}
%Let $s_1 < s_2 < \ldots$ be the times at which $v$ fires, and $t_1 < t_2 < \ldots$ the times at which $w$ fires.  Then we must show
%	\[ t_1 \leq s_1 \leq t_2 \leq s_2 \leq \ldots. \]
Induct on $t$.  If $u_t(\sigma,w) = u_t(\sigma,v)$, then by (\ref{laplacianofodometer})
	\[ U^t\sigma(v) - U^t\sigma(w) = \sigma(v) - \sigma(w) \leq 0. \]
Thus if $v$ fires at time $t$, then $w$ must also fire, hence 
%(\ref{interlacedodometers}) holds at time $t+1$.  
	\begin{equation} \label{interlacedodometers} u_{t+1}(\sigma,v) \leq u_{t+1}(\sigma,w) \leq u_{t+1}(\sigma,v)+1. \end{equation}
On the other hand, if $u_t(\sigma,w) = u_t(\sigma,v)+1$, then since $\sigma$ is confined, we have by (\ref{laplacianofodometer})
	\[ U^t\sigma(v) - U^t\sigma(w) = n + \sigma(v) - \sigma(w) \geq 0. \]
Thus if $w$ fires at time $t$, then $v$ must also fire, so once again (\ref{interlacedodometers}) holds.
%By Lemma~\ref{whichsitesareunstable},
%	\[ u_t(\sigma,v) = \#\{s<t | \sigma(v) \equiv k (\mod n) \mbox{ for some } \alpha_s < k \leq \alpha_{s+1} \}. \]
\end{proof}

\begin{lemma}
\label{divisibilitycriterion}
If $\sigma$ is confined, then $U^t \sigma = \sigma$ if and only if $n | \alpha_t$.
\end{lemma}

\begin{proof}
If $U^t \sigma=\sigma$, then $n|\alpha_t$ by Lemma~\ref{constantodometer}.  For the converse, 
write $\ell = \min_v u_t(\sigma,v)$.  By Lemma~\ref{interlacing} we have 
	\[ u_t(\sigma,v) \in \{\ell, \ell+1\} \] 
for all vertices $v$.  If $n | \alpha_t$, then since
	\[ \alpha_t = \sum_{v=1}^n u_t(\sigma,v), \]
we must have $u_t(\sigma,v) = \ell$ for all $v$, so $U^t \sigma =\sigma$ Lemma~\ref{constantodometer}.
\end{proof}

Let $\sigma$ be a confined state on $K_n$.  By the pigeonhole principle, there exist times $0 \leq s<t \leq n$ with
	\[ \alpha_s \equiv \alpha_t \modulo{n}. \]
By Lemma~\ref{divisibilitycriterion} it follows that $U^s \sigma = U^t \sigma$, so~$\sigma$ has eventual period at most~$n$.  

Our next result improves this bound a bit.  Write $m(\sigma)$ for the eventual period of $\sigma$, and 
	\[ \nu(\sigma) = \#\{ \sigma(v) | v\in [n] \} \]
for the number of distinct heights in $\sigma$.

\begin{prop}
\label{distinctheights}
For any chip configuration $\sigma$ on $K_n$,
	\[ m(\sigma) \leq \nu(\sigma). \]
\end{prop}

\begin{proof}
If $a(\sigma)=1$, then $m(\sigma)=1$.
Since $\nu(U\sigma) \leq \nu(\sigma)$ and $m(U\sigma) = m(\sigma)$, by Lemma~\ref{eventuallyconfined} it suffices to prove the proposition for confined states~$\sigma$.
Moreover, by relabeling the vertices we may assume that $\sigma(1) \geq \sigma(2) \geq \ldots \geq \sigma(n)$. 
%Let $v_1=1$ and
%	\[ v_j = \min \{v>v_{j-1}| \sigma(v)<\sigma(v-1)\}, \qquad j=2, \ldots, \nu. \]
%Relabel the vertices so that $\sigma(v) \geq \sigma(w)$ for $v\leq w$.  
Write the total activity $\alpha_t$ as
	\[ \alpha_t = Q_t n + R_t \]
for nonnegative integers $Q_t$ and $R_t$ with $R_t \leq n-1$.  By Lemma~\ref{interlacing}, since $\sum_v u_t(\sigma,v) = \alpha_t$, we must have
	\[ u_t(\sigma,v) = \begin{cases} Q_t+1, & v \leq R_t \\
							Q_t, & v > R_t.  \end{cases} \]
Moreover if $R_t>0$, then $\sigma(R_t)>\sigma(R_t+1)$.  Since the set
	\[ V = \{0\} \cup \{ v\in [n-1] | \sigma(v) > \sigma(v+1) \} \]
has cardinality $\nu(\sigma)$, by the pigeonhole principle there exist times $0 \leq s <t \leq \nu(\sigma)$ with $R_s=R_t$.  Then
	\[ u_t(\sigma,v) - u_s(\sigma,v) = Q_t - Q_s \]
independent of $v$, hence $U^s \sigma = U^t \sigma$ by Lemma~\ref{constantodometer}.  Thus $m(\sigma) \leq t-s \leq \nu(\sigma)$.
\end{proof}

Bitar \cite{Bitar} conjectured that any parallel chip-firing configuration on a connected graph of~$n$ vertices has eventual period at most~$n$.  A counterexample was found by Kiwi et al.\ \cite{KNTG}.  Anne Fey (personal communication) has found a counterexample on a regular graph.  Proposition~\ref{distinctheights} shows that Bitar's conjecture holds on the complete graph; it is also known to hold on trees \cite{BG} and on cycles \cite{Dall'Asta}.
It would be interesting to find other classes of graphs for which  Bitar's conjecture holds.

Next we relate the period of a chip configuration to its activity.  We will use the fact that the rotation number of a circle map determines the periods of its periodic points: if $f : \R \to \R$ is a monotone degree one lift (\ref{degreeonelift}) with $\rho(f)=p/q$ in lowest terms, then all periodic points of $\bar{f} : S^1 \to S^1$ have period $q$; see Proposition~4.3.8 and Exercise~4.3.5 of \cite{HK}.

Note that if a chip configuration $\sigma$ on $K_n$ is periodic of period $m$, and $a(\sigma)=p/q$, then
	\begin{equation} \label{alpham} \frac{\alpha_m}{mn} = \frac{p}{q} \end{equation}
and $\alpha_{km} = k\alpha_m$ for all $k \in \N$.

\begin{lemma}
\label{denominatorofactivity}
If $a(\sigma) = p/q$ and $(p,q)=1$, then $m(\sigma)=q$.
\end{lemma}

\begin{proof}
Since $a(U\sigma)=a(\sigma)$ and $m(U\sigma)=m(\sigma)$, we may assume that~$\sigma$ is periodic.  Now if $a(\sigma)=1$, then $m(\sigma)=1$, so we may assume $a(\sigma)<1$.  In particular, $\sigma$ is confined by Lemma~\ref{eventuallyconfined}.  Write $m=m(\sigma)$.  By Lemma~\ref{constantodometer} and (\ref{alpham})
	\[ u_m(\sigma,v) = \frac{\alpha_m}{n} = \frac{pm}{q} \]
for all vertices $v$.  Hence $q|pm$ and hence $q|m$.

Let $f$ be the lift of the circle map $\bar{f}$ associated to $\sigma$.  By Lemma~\ref{activityequalsrotationnumber} we have $\rot(f) = a(\sigma)$, so all periodic points of $\bar{f}$ have period~$q$.  By Lemma~\ref{theiterates} we have
	\[ f^{qm}(0) = \frac{\alpha_{qm}}{n} = \frac{q\alpha_m}{n} = pm \in \Z, \]
so $0$ is a periodic point of $\bar{f}$.  Thus $f^q(0) \in \Z$.  Now by Lemma~\ref{theiterates}
	\[  \alpha_q = n f^q(0) \]
hence $U^q \sigma = \sigma$ by Lemma~\ref{divisibilitycriterion}.  Thus $m|q$.
\end{proof}

Given $1 \leq p<q \leq n$ with $(p,q)=1$ and $p/q \leq 1/2$, one can check that the chip configuration $\sigma$ on $K_n$ given by
	\[ \sigma(v) = \begin{cases} v+p-1, & v \leq q-1-p \\
						   v+n+p-q-1, & q-p \leq v \leq q-1 \\
						   n+p-1, & v \geq q. \end{cases}  \]
%
% possible alternatives:
%
% does this work?
%	\[ \sigma(v) = \begin{cases} v + n + p- q -1 , & v=1,\ldots, q \\
%						   n+2p-q, & v=q+1, \ldots n. \end{cases} \]
%
% this one doesn't:
%	\[ \sigma(v) = v + n - q + p -1 , \qquad v=1,\ldots, q \]
% and $\sigma(v) =n-1$ for $v \geq q+1$.
%   						   
has activity $a(\sigma)=p/q$.  For a similar construction on more general graphs in the case $p=1$, see \cite[Prop.\ 6.8]{Prisner}.  By Lemma~\ref{denominatorofactivity}, $m(\sigma)=q$.
So for every integer $q=1,\ldots,n$ there exists a chip configuration on $K_n$ of period $q$.  

Despite the existence of states of period as large as $n$, states of smaller period are in some sense more prevalent.  One way to capture this is the following.

\begin{theorem}
\label{lotsofsmallperiods}
If $\sigma_2,\sigma_3,\ldots$ is a sequence of chip configurations satisfying (\ref{stability}), (\ref{converginghistograms}) and (\ref{notidentity}), then
for each $q\in \N$ there is a constant $c=c_q >0$ such that for all sufficiently large $n$, at least $c n$ of the states $\{\sigma_n+k\}_{k=0}^n$ have eventual period $q$.
\end{theorem}

\begin{proof}
Fix $\epsilon>0$, and choose $p<q$ with $(p,q)=1$.  By Proposition~\ref{devil'sstaircase} we have $s^{-1}(p/q) = [a_p,b_p]$ for some $a_p<b_p$.  By Proposition~\ref{strongerconvergence}, for sufficiently large $n$ we have 
	\[ a(\sigma_n + \floor{ny}) = s_n(y) = p/q \]
for all $y \in [a_p+\epsilon/2q,b_p-\epsilon/2q]$.  Thus all states $\sigma_n+k$ with
	\[ \left(a_p+\frac{\epsilon}{2q}\right) n < k < \left(b_p-\frac{\epsilon}{2q}\right) n \]
have activity $p/q$ and hence, by Lemma~\ref{denominatorofactivity}, eventual period $q$.  We can therefore take
	\[ c_q = \sum_{p : (p,q)=1} (b_p - a_p)  - \epsilon \]
for any $\epsilon>0$.
\end{proof}

The rest of this section is devoted to proving the existence of a \emph{period~$2$ window}: any chip configuration on $K_n$ with total number of chips strictly between $n^2-n$ and $n^2$ has eventual period~$2$. 

The following lemma is a special case of \cite[Prop.\ 6.2]{Prisner}; we include a proof for completeness.  Write $r(\sigma)$ for the number of unstable vertices in $\sigma$.  

\begin{lemma}
\label{reflectionsymmetry}
If $\sigma$ and $\tau$ are chip configurations on $K_n$ with $\sigma(v) + \tau(v) = 2n-1$ for all $v$, then $a(\sigma)+a(\tau)=1$.
\end{lemma}

\begin{proof}
If $\sigma(v) + \tau(v) = 2n-1$ for all $v$, then for each vertex $v$, exactly one of $\sigma(v), \tau(v)$ is $\geq n$.  Hence $r(\sigma)+r(\tau)=n$, and
	\begin{align*} U\sigma(v) + U\tau(v) &= \sigma(v)+r(\sigma) + \tau(v) + r(\tau) - n \\ &= 2n-1 \end{align*}
for all vertices $v$.  Inducting on $t$, we obtain
	\[ r(U^t \sigma) + r(U^t \tau) = n \]
for all $t \geq 0$, and hence $a(\sigma)+a(\tau)=1$ by (\ref{theactivity}).
%	\[ a(\sigma) + a(\tau) = \lim_{t \to \infty} \frac{1}{nt} \sum_{s=0}^{t-1} r(U^s \sigma) + r(U^s \tau) = 1. \qed \]
%\renewcommand{\qedsymbol}{}
\end{proof}

%If $\sigma(1) \geq \sigma(2) \geq \ldots \geq \sigma(n)$, write
Given a chip configuration $\sigma$ on $K_n$, for $j=1,\ldots,n$ we define \emph{conjugate} configurations
	\[ c^j \sigma (v) = \begin{cases} \sigma(v)+j-n, & v \leq j \\
							\sigma(v)+j, & v>j. \end{cases} \]
Note that $c^n \sigma = \sigma$.
%\begin{lemma}
%	\[\alpha_{t,i} = g(\alpha_{t-1,i}+i)-i. \]
%\end{lemma}
%
%\[ f^t \left(\frac{i}{n}\right) = \frac{\alpha_{t,i}+i}{n} \]
%
%\[ a(c^i \sigma) = \lim_{t \to \infty} \frac{\alpha_{t,i}}{nt} = \lim_{t \to \infty} \frac{f^t(i/n)}{t} = \rho(f) = a(\sigma). \]

\begin{lemma}
Let $\sigma$ be a chip configuration on $K_n$, and fix $j \in [n]$.
For all~$t\geq 0$, we have for $v \leq j$
	\[ u_t(\sigma, v)-1 \leq u_t(c^j \sigma,v) \leq u_t(\sigma,v), \] 
while for $v > j$
	\[ u_t(\sigma, v) \leq u_t(c^j \sigma,v) \leq u_t(\sigma,v)+1. \] 
\end{lemma}

\begin{proof}
%maybe can reduce to case j=1?  not sure if this really simplifies things much.
Induct on $t$.   By (\ref{laplacianofodometer})
	\begin{equation} \label{morelaplacian} U^t c^j \sigma(v) = c^j \sigma(v) + \alpha_{t,j} - n\, u_t(c^j\sigma,v) \end{equation}
where
	\[ \alpha_{t,j} = \sum_{w=1}^n u_t(c^j\sigma,w). \]
By the inductive hypothesis,
	\begin{equation} \label{activitybounds} -j \leq \alpha_{t,j} - \alpha_{t,0} \leq n-j. \end{equation}

Fix a vertex $v$, and write 	
	\[ b = u_t(c^j\sigma,v) - u_t(\sigma,v). \]  
By the inductive hypothesis we have $b \in \{-1,0,1\}$.  From (\ref{morelaplacian}) we have
	\begin{align} \label{heightdifference}
	U^t c^j \sigma (v) - U^t \sigma(v) &= \alpha_{t,j} - \alpha_{t,0} + c^j \sigma(v) - \sigma(v) - n b 
%	\nonumber \\
%		&= \alpha_{t,j} - \alpha_{t,0} + j - nb'
	 \end{align}
%where
%	\[ b' = \begin{cases} b+1, & v\leq j \\ b, & v>j.  \end{cases} \]

Case 1:  $v\leq j$ and $b=-1$.  By (\ref{activitybounds}), the right side of (\ref{heightdifference}) is $\geq 0$, so if $v$ is unstable in $U^t \sigma$, then it is also unstable in $U^t c^j \sigma$.  Thus
	\begin{equation} \label{smallerodometer} u_{t+1}(\sigma, v)-1 \leq u_{t+1}(c^j \sigma,v) \leq u_{t+1}(\sigma,v), \end{equation}
completing the inductive step.

Case 2: $v\leq j$ and $b=0$.  By (\ref{activitybounds}), the right side of (\ref{heightdifference}) is $\leq 0$, so if~$v$ is unstable in $U^t c^j \sigma$, then it is also unstable in $U^t \sigma$.  Thus once again (\ref{smallerodometer}) holds.
	
Case 3: $v>j$ and $b=0$.  By (\ref{activitybounds}), the right side of (\ref{heightdifference}) is $\geq 0$, so 
		\begin{equation} \label{biggerodometer} u_{t+1}(\sigma, v) \leq u_{t+1}(c^j \sigma,v) \leq u_{t+1}(\sigma,v)+1, \end{equation}
completing the inductive step.  

Case 4: $v>j$ and $b=1$.  By (\ref{activitybounds}), the right side of (\ref{heightdifference}) is $\leq 0$, so once again (\ref{biggerodometer}) holds.
\end{proof}

\begin{corollary}
\label{conjugateactivities}
For any chip configuration $\sigma$ on $K_n$ and any $j \in [n]$,
 \[ a(c^j \sigma) = a(\sigma). \]
\end{corollary}

It turns out that the circle maps corresponding to $\sigma$ and $c^j \sigma$ are conjugate to one another by a rotation.  This gives an alternative proof of the corollary, in the case when both $\sigma$ and $c^j \sigma$ are preconfined.

\begin{lemma}
\label{ican'tbelievethisdoesn'tfollowfromsomethingwealreadyproved}
Let $\sigma$ be a chip configuration on $K_n$.
If $u_2(\sigma,v) \geq 1$ for all $v$, then 
	$u_{2t}(\sigma,v) \geq t$ 
for all $v$ and all $t \geq 1$.
\end{lemma}

\begin{proof}
Induct on $t$.  By the inductive hypothesis, $\alpha_{2t} \geq nt$.  Fix a vertex~$v$, and suppose that $u_{2t}(\sigma,v)=t$.  Since $u_2(\sigma,v)\geq 1$, either $\sigma(v) \geq n$ or $U\sigma(v) \geq n$.  In the former case, by (\ref{laplacianofodometer})
	\[ U^{2t} \sigma(v) - \sigma(v) = \alpha_{2t} - nt \geq 0, \]
so $v$ fires at time $2t$.  Hence $u_{2t+1}(\sigma,v) \geq t+1$.  Summing over $v$ yields
	\[ \alpha_{2t+1} \geq \alpha_1 (t+1) + (n-\alpha_1) t = nt + \alpha_1. \]

  In the latter case, if $v$ does not fire at time $2t$, we have by (\ref{laplacianofodometer}) 
	\[ U^{2t+1} \sigma(v) - U\sigma(v) = \alpha_{2t+1} - \alpha_1 - nt \geq 0, \]
so $v$ fires at time $2t+1$.  Hence $u_{2t+2}(\sigma,v) \geq t+1$ for all $v$.
\end{proof}

Given a chip configuration $\sigma$ on $K_n$, write
	\[ \ell(\sigma) = \min \{ \sigma(1),\ldots,\sigma(n) \} \]
and 
	\[ r(\sigma) = \# \{v\in [n] : \sigma(v) \geq n\}. \]
Write
	\[ |\sigma| = \sum_{v=1}^n \sigma(v) \]
for the total number of chips in the system.

\begin{theorem}
\label{period2window}
Every chip configuration $\sigma$ on $K_n$ with $n^2 - n < |\sigma| < n^2$ has eventual period~$2$.
\end{theorem}

\begin{proof}
If $|\sigma |<n^2$, then $a(\sigma)<1$, so by Lemma~\ref{eventuallyconfined} we may assume~$\sigma$ is confined. 
% NOTE c^j \sigma may not be confined! 
Moreover, by relabeling the vertices we may assume that $\sigma(1) \geq \sigma(2) \geq \ldots \geq \sigma(n)$.  In particular, for any vertex $v$, 
% not needed:
%
%since $\sigma(w) \leq \sigma(v)+n-1$ for $w=1,\ldots,v-1$, we have
%	\begin{align*} |\sigma| &\leq (v-1)(\sigma(v)+n-1) + (n-v+1)\sigma(v) \\
%			&= n\sigma(v) + (v-1)(n-1). \end{align*}
%Thus if $|\sigma|>n^2-n$, we obtain
%	\begin{equation} \label{lowerheightbound} \sigma(v) > \frac{(n-v+1)(n-1)}{n} \geq n-v. \end{equation}
%Likewise, 
since $\sigma(w) \geq \sigma(v) - n +1$ for $w=v+1,\ldots,n$, we have
	\begin{align*} |\sigma| &\geq v\sigma(v) + (n-v)(\sigma(v)-n+1) \\
			&= n\sigma(v) - (n-v)(n-1). \end{align*}
Since $|\sigma|<n^2$, we obtain
	\begin{equation} \label{upperheightbound} \sigma(v) < \frac{n^2+(n-v)(n-1)}{n} \leq 2n-v. \end{equation}

For a fixed vertex $v$, as $j$ ranges over $[n]$, the quantity $c^j \sigma (v)$ takes on each of the values
	\[ \sigma(v)+v-n, \ldots, \sigma(v)+v-1 \] 
exactly once.  By (\ref{upperheightbound}), at least $\sigma(v)+v-n$ of these values are $\geq n$, hence
	\begin{equation} \label{sumofrs} \sum_{j=1}^n r(c^j \sigma) \geq \sum_{v=1}^n (\sigma(v)+v-n) = |\sigma| - \frac{n(n-1)}{2}. \end{equation}
Since $\ell(c^j \sigma) = \sigma(j)+j-n$, we have
	\begin{equation} \label{sumofells} \sum_{j=1}^n \ell(c^j \sigma) = |\sigma| - \frac{n(n-1)}{2}. \end{equation}

Summing equations (\ref{sumofrs}) and  (\ref{sumofells}), and using $|\sigma|>n^2-n$, we obtain
	\[ \sum_{j=1}^n (\ell(c^j \sigma) + r(c^j \sigma)) > n^2 -n. \]
Since each term in the sum on the left is a nonnegative integer, we must have
	\[ \ell(c^j \sigma) + r(c^j \sigma) \geq n. \]
for some $j \in [n]$.  Thus every vertex $v$ fires at least once during the first two updates of $c^j\sigma$.  From Corollary~\ref{conjugateactivities} and Lemma~\ref{ican'tbelievethisdoesn'tfollowfromsomethingwealreadyproved} we obtain
	\[ a(\sigma) = a(c^j \sigma) \geq \frac12. \]

Finally, the chip configuration $\tau(v) := 2n-1-\sigma(v)$ also satisfies $n^2 - n < |\tau| < n^2$, so $a(\tau) \geq \frac12$.  By Lemma~\ref{reflectionsymmetry} we have $a(\sigma) + a(\tau)=1$, so $a(\sigma) = a(\tau) = \frac12$.  From Lemma~\ref{denominatorofactivity} we conclude that $m(\sigma)=2$.
\end{proof}

\section*{Acknowledgements}
The author thanks Anne Fey for many helpful conversations, and an anonymous referee for suggesting the formulation of generalized chip configurations in terms of measures on $[0,2)$.

\end{document}